\documentclass[twoside,10pt]{amsart}
\usepackage{amsmath}%
\usepackage{amsfonts}%
\usepackage{amssymb}%
\usepackage{graphicx}%

\newtheorem{theorem}{Theorem}

\newtheorem{corollary}[theorem]{Corollary}

\newtheorem{definition}[theorem]{Definition}
\newtheorem{example}[theorem]{Example}
\newtheorem{examples}[theorem]{Examples}

\newtheorem{lemma}[theorem]{Lemma}

\newtheorem{proposition}[theorem]{Proposition}
\newtheorem{remark}[theorem]{Remark}

\newenvironment{Rmk}{\begin{remark}\em}{\end{remark}}
\newenvironment{Cor}{\begin{corollary}\em}{\end{corollary}}

\newtheorem{prf}{Proof}

\newenvironment{Prf}{\begin{prf}\em}{\qed\end{prf}}
\newtheorem{prff}{}

\newcommand{\NOT}[1]{}
\newcommand{\pa}{\par\medskip}

\newcommand{\RRe}{\mathrm{Re}}

\newcommand{\eps}{\varepsilon}

\newcommand{\U}{\mathbf{1}}

\begin{document}

%Topmatter

\title[Operators between Hilbert Spaces Viewed as Only Linear Topological]{Operators between Hilbert Spaces Viewed as Only Linear Topological -- Towards a Classification}
%\author{Eliahu Levy, Mathematics, Technion, Haifa}

%Author info
\author{Eliahu Levy}
\address{Department of Mathematics,
Technion -- Israel Institute of Technology,
Haifa 32000, Israel}
\email{eliahu@math.technion.ac.il}
\date{}


%\thanks{}
%\keywords{}
%\subjclass{}

%\date{ \ \\
%\date{\dateline{2008}{2008}\\
%\small Mathematics Subject Classification: }

%End topmatter

\begin{abstract}
In \textit{topological equivalence}, a bounded linear operator between Banach spaces -- we focus on the case of Hilbert spaces -- is viewed as only acting linearly and continuously between them \textit{qua} different spaces with the structure of linear topological space. For instance, invertible operators in Banach spaces (that is, isomorphisms among them) will make up one equivalence class for each class of isomorphic spaces. On the other hand, compact and non-compact operators, or operators with or without a kernel, clearly will not. We make some crucial steps towards describing invariants that will characterize these topological equivalence classes
\end{abstract}

\maketitle

\tableofcontents

\section{Introduction, Topological Equivalence of Operators}
In \textit{topological equivalence} of operators between Banach spaces -- we focus on the case of Hilbert spaces -- we wish to view an operator as only acting linearly and continuously between the spaces, \textit{qua} \textit{different} spaces with the structure of linear topological space. Thus
\begin{definition}
Two operators $T_1:X_1\to Y_1$ and $T_2:X_2\to Y_2$ between Banach spaces will be called \textit{Topologically equivalent}
if there exist isomorphisms onto $U:X_1\to X_2$ and $V:Y_1\to Y_2$ so that $T_2=V^{-1}T_1U$.
\end{definition}
Note that even if $X_1=X_2=Y_1=Y_2$, $U$ and $V$ are allowed to be \textit{different} isomorphisms.\pa

Of course, for $T_1$ and $T_2$ to be topologically equivalent, $X_1$ and $X_2$, and likewise $Y_1$ and $Y_2$,
must be isomorphic.\pa

We make some crucial steps towards describing invariants that will characterize the topological equivalence classes.\pa

First, any operator $T:\mathcal{H}\to\mathcal{H}'$ between Hilbert is topologically (even `isometrically') equivalent to $A:=\sqrt{T^\star T}$ -- a positive semi-definite operator in $\mathcal{H}$. Indeed $\|Tx\|=\|Ax\|,\,x\in\mathcal{H}$ ,\pa

One focuses on $I$-subspaces, where $I$ is the interval in the positive half-line, consisting of the numbers $<$,\,$\le$,\,$>$, or $\ge$ from some fixed $\mu$. An $I$-subspace $Y$ for $T$ (here the same for $A:=\sqrt{T^\star T}$ mentioned above) is then defined by requiring that $\|Tx\|/\|x\|\in I$ for all $0\neq x\in Y$.\pa

And we prove that for maximal $[0,\mu]$-subspaces, and for $\mu>\mu'$ the quotient of a maximal $[0,\mu]$-subspaces over some maximal $[0,\mu']$-subspace contained in it, are all isomorphic, a special instant featuring in the spectral decomposition of the a positive semi-definite $A$.\pa

That allows us to relate features of that spectral decomposition to notions more unaffected when a Hilbert norm is replaced by an equivalent Hilbert norm, thus to formulate some necessary condition for topological equivalence, Theorem
\ref{th:dm-eq}.\pa

To obtain sufficient conditions, one notes that an operator is characterized up to topological equivalence (indeed
isometric) by the dimensions of the parts in the latter spectral decomposition -- these forming a cardinal number-valued ($\sigma$-additive) `measure' on the Borel subsets of $\mathbb{R}^+=[0,\infty[$.\pa

And one sees what manipulations in the latter preserve topological equivalence, toward attaching to an operator
what we call an \textit{eigenvalue sequence with infinity intermissions} with some criterions for topological equivalence, Theorem \ref{th:cond}.\pa

\section{$I$-Spaces for Operators, Reduction to Positive Semi-Definite}\label{s:A}
Let $T:X_1\to X_2$ be a continuous linear operator between Banach spaces.\pa
\begin{definition}
Let $I$ be an interval of the form $[0,\mu]$ or $[0,\mu[$, $\mu\ge0$, in $\mathbb{R}^+=[0,+\infty[$, or of the `complement' form $]\mu,\infty[$ or $[\mu,\infty[$. We say that a linear subspace $Y\subset X_1$ (here, not necessarily closed) is an \textit{$I$-subspace (with respect to $T$)} if for every $0\ne x\in Y$\,\,$\|Tx\|/\|x\|\in I$.
\end{definition}

In particular, one may consider \textit{maximal} (by inclusion) {$I$-subspaces} in $X$ with respect to $T$ .\pa

Clearly, if $I$ is closed, (i.e.\ is of the form $[0,\mu]$ or $[\mu,\infty[$) then the closure of an $I$-subspace is an $I$-subspace. Therefore \textit{if $I$ is closed a maximal $I$-subspace is closed}.\pa

By a straightforward application of Zorn's lemma, \textit{every $I$-subspace is contained in a maximal $I$-subspace}.\pa

Let now $\mathcal{H}$ and $\mathcal{H}'$ be Hilbert spaces (assumed complex, or if real complexify them, since we shall need the spectral decomposition) and let $T:\mathcal{H}\to\mathcal{H}'$ be a bounded linear
operator.\pa

We have for any $x\in\mathcal{H}$
$$\|Tx\|^2=\langle Tx,Tx\rangle=\langle T^\star Tx,x\rangle,$$
where $T^\star T$ is a positive semi-definite operator from $\mathcal{H}$ to itself.\pa

$T^\star T$ has a unique positive semi-definite square root $A=\sqrt{T^\star T}$, and we have, for any $x\in\mathcal{H}$
$$\|Tx\|^2=\langle T^\star Tx,x\rangle=\langle A^2x,x\rangle=\langle Ax,Ax\rangle=\|Ax\|^2.$$
(As a consequence $\mathrm{Ker}T=\mathrm{Ker}A$).\pa

We shall thus have an isometry $U$ between the image of $A$ in $\mathcal{H}$ and the image of $T$ in
$\mathcal{H}'$, by mapping $Ax\mapsto Tx$. This isometry can, of course, be uniquely extended to an isometry
$U$ from the \textit{closure} $\mathrm{cl}(A\,\mathcal{H})$ in $\mathcal{H}$ to the closure $\mathrm{cl}(T\,\mathcal{H})$ in $\mathcal{H}'$. Of course $T=UA$.\pa

Thus, $T$, as operator from $\mathcal{H}\to\mathrm{cl}(T\,\mathcal{H})$, is \textit{topologically equivalent} to
$A=\sqrt{T^\star T}$ from $\mathcal{H}\to\mathrm{cl}(A\mathcal{H})$.\pa

Hence, when the linear topological spaces can be normed as Hilbert, we, up to topological equivalence, may restrict ourself to \textit{positive semi-definite operators in a Hilbert space $\mathcal{H}$ as operators from $\mathcal{H}$ to the closure of their image}.\pa

A particular case is when $T$, consequently $A^2$ and $A$, have no kernel. Since $A$ is Hermitian, for it \textit{kernel $=$ orthogonal complement of the image}, so in that case it has a dense image in $\mathcal{H}$. Thus $U$ is an isometry from $\mathcal{H}$ \textit{onto} $\mathcal{H}'$. Consequently, in this case ($\mathcal{H}\simeq\mathcal{H}'$) and $T$ is topologically equivalent to $A$.\pa

Returning to the general case, $A$, a (semi-definite) Hermitian operator, has a \textit{spectral decomposition} \cite{DunfordSchwartz}, given by a `decomposition of unity', i.e.\ a spectral measure on $\mathcal{H}$,\,\,$d\mathcal{E}(\lambda)$, defined on $\mathbb{R}^+=[0,+\infty[$,%
\footnote{That is, $\mathcal{E}(\tau)$ is an (operator on $\mathcal{H}$)-valued measure on the Borel subsets $\tau$
of $\mathbb{R}^+$, strongly $\sigma$-additive, whose value on sets is always an orthogonal projection, and such that
the measure of the empty set is the operator $0$, the measure of the whole $\mathbb{R}^+$ is the operator $\U$, and
the measure of the intersection of two Borel sets is the product of the measures of the sets, which thus commute.
In particular, the measures of disjoint sets are orthogonal projections on orthogonal subspaces.}
such that $A=\int_{[0,+\infty[}\lambda\,d\mathcal{E}(\lambda)$, i.e., $\forall\,x,x'\in\mathcal{H}$
$$\langle Ax,x'\rangle=\int_{[0,+\infty[}\lambda\langle d\mathcal{E}(\lambda)x,x'\rangle,$$
which implies, by the rules of manipulating spectral measures, that
$A^2=\int_{[0,+\infty[}\lambda^2\,d\mathcal{E}(\lambda)$, i.e.,
$$\langle T^\star Tx,x'\rangle=\langle A^2x,x'\rangle=\int_{[0,+\infty[}\lambda^2\langle d\mathcal{E}(\lambda)x,x'\rangle.$$

\section{Isomorphism Property of Maximal $[0,\mu]$-Spaces and their Quotients}\label{s:iso}
\begin{proposition}\label{pr:Ispc}
Referring to the above spectral decomposition of $A=\sqrt{T^\star T}$, let $I$ be of the form $[0,\mu]$, $[0,\mu[$, $]\mu,\infty[$ or $[\mu,\infty[$, $\mu\ge0$.

Then $\mathcal{E}(I)\mathcal{H}$ is a maximal $I$-subspace with respect to $T$, equivalently with respect to $A$.\pa
\end{proposition}
\begin{Prf}
Let $x\in\mathcal{E}(I)\mathcal{H}$. Then $x$ is orthogonal to $\mathcal{E}(I^c)\mathcal{H}$ and we have (here $I^2$ denotes the set of squares of members of $I$),
\begin{eqnarray*}
&&\|Tx\|^2=\langle A^2x,x\rangle=\int_{[0,+\infty[}\lambda^2\langle d\mathcal{E}(\lambda)x,x\rangle=
\int_I\lambda^2\langle d\mathcal{E}(\lambda)x,x\rangle\in\\
&&\in I^2\cdot\int_I\langle d\mathcal{E}(\lambda)x,x\rangle=I^2\cdot\langle\mathcal{E}(I)x,x\rangle=I^2\cdot\|x\|^2.
\end{eqnarray*}

Consequently $\|Tx\|\in I\cdot\|x\|$, making $\mathcal{E}(I)\mathcal{H}$ an $I$-subspace.\pa

Now denote by $I^c$ the complement of $I$ in $\mathbb{R}^+=[0,\infty[$. Since $I^c$ is also an interval of our form, $\mathcal{E}(I^c)\mathcal{H}$ is an $I^c$-subspace  with respect to $A=\sqrt{T^\star T}$, thus to $T$.\pa

And we have the orthogonal decomposition $\mathcal{H}=\mathcal{E}(I)\mathcal{H}\oplus\mathcal{E}(I^c)\mathcal{H}$.\pa

Now suppose $\mathcal{E}(I)\mathcal{H}$ were not maximal, i.e.\ contained in a larger $I$-subspace $Y$. Then $Y$ must intersect $\mathcal{E}(I^c)\mathcal{H}$ nontrivially, i.e.\ there would be a nonzero
$a\in Y\cap\mathcal{E}(I^c)\mathcal{H}$. Then $\|Ta\|/\|a\|$ must be both in $I$ and $I^c$ -- impossible.
\end{Prf}

\begin{Rmk}
For our record here, note trivially that for $\mu=0$, $I=[0,0]=\{0\}$, and we have: for $Y$ to be an $I$-subspace means that $T$ vanishes on $Y$; for $Y$ to be a \textit{maximal} $I$-subspace means that $Y=\mathrm{Ker}\,T$; and the latter is $\mathcal{E}(\{0\})\mathcal{H}$ in the above spectral decomposition.
\end{Rmk}

\begin{proposition}
Suppose $I$ is an interval of the form $[0,\mu]$,\,$\mu\ge0$.

($I$ being closed, clearly the closure $\mathrm{cl}\,Y$ of an $I$-subspace $Y$ (with respect to $T$, equivalently to $A=\sqrt{T^\star T}$) is also an $I$-subspace. Therefore \textit{a maximal $I$-subspece must be closed}.)

\textbf{(i)} Let $Y$ be an $I$-subspace.

Then the restriction to $Y$ of the orthogonal projection on $\mathcal{E}(I)\mathcal{H}$ is \textit{injective} on $Y$.

\textbf{(ii)} Let $Y$ be a \textit{maximal} $I$-subspace (hence closed).

Then the restriction to $Y$ of the orthogonal projection on $\mathcal{E}(I)$ is also onto, hence is an \textit{isomorphism} between $Y$ and $\mathcal{E}(I)\mathcal{H}$.
\end{proposition}
\begin{Prf}

\textbf{(i)} Recall that we denote by $I^c$ the complement of $I$ in $\mathbb{R}^+=[0,\infty[$, here
$I^c=]\mu,\infty[$.

We have the orthogonal decomposition $\mathcal{H}=\mathcal{E}(I)\mathcal{H}\oplus\mathcal{E}(I^c)\mathcal{H}$, where, by Proposition \ref{pr:Ispc}, $\mathcal{E}(I)\mathcal{H}$ is an $I$-subspace and $\mathcal{E}(I^c)\mathcal{H}$ is an $I^c$-subspace with respect to $T$, equivalently with respect to $A=\sqrt{T^\star T}$. Therefore, if $Y$ and $\mathcal{E}(I^c)\mathcal{H}$ had a nontrivial intersection, i.e.\ there was an $a\ne0$ contained in both, then $\|Ta\|/\|a\|$ would have been both in $I$ and $I^c$ -- impossible.

Consequently $Y\cap\mathcal{E}(I^c)=0$, making the restriction to $Y$ of the orthogonal projection on $\mathcal{E}(I)$ in the above orthogonal decomposition injective.

\medskip

\textbf{(ii)} Suppose now that $Y$ is a maximal $I$-subspace, and, toward a contradiction, that the projection is not onto, i.e.\ there is an $a\in\mathcal{E}(I)$ not in the image of the above projection.\pa

Let $y$ be any member of $Y$.

Denote by $y_{I}$ and $y_{I^c}$ the projections of $y$ on $\mathcal{E}(I)\mathcal{H}$ and $\mathcal{E}(I^c)\mathcal{H}$, resp.

Now, $a+y_{I}$ is in $\mathcal{E}(I)$ and $y_{I^c}$ is in $\mathcal{E}(I^c)$, which are orthogonal and
$A$-invariant. Therefore the above projection sends $a+y$ to $a+y_{I}$, hence different from all projections of
members of $Y$, implying $a\notin Y$.\pa

Thus, the linear space $Y'$ generated by $Y$ and $a$ is strictly bigger than $Y$. Therefore, as $y$ is a maximal
$I$-subspace, $Y'$ must contain some $a+y$ with $\|A(a+y)\|\notin I\|a+y\|$\,\,\,(while $\|A(a)\|\in I\|a\|$ since
$a\in\mathcal{E}(I)\mathcal{H}$, and $\|A(y)\|\in I\|y\|$ since $y\in Y$).\pa

So (recall that $I=[0,\mu]$) we have\pa
$$\|A(a+y)\|^2>\mu^2\|a+y\|^2,$$
i.e.,
\begin{equation}\label{eq:ay}
\mu^2\|a\|^2-\|Aa\|^2+\mu^2\|y\|^2-\|Ay\|^2+
2\RRe\Big(\mu^2\langle a,y\rangle-\langle Aa,Ay\rangle\Big)<0.
\end{equation}

On the other hand, consider the segment $\{ta+(1-t)y\,|\,0\le t\le 1\}$. The quadratic expression in $t$
$$\|A(ta+(1-t)y)\|^2-\mu^2\|ta+(1-t)y\|^2$$
is (recall that $I=[0,\mu]$),\quad $\le0$ for $t=0$ (since $a\in\mathcal{E}(I)\mathcal{H}$) and for $t=1$ (since $y\in Y$), but positive for $t=1/2$. Therefore its derivative at $t=0$ must be positive.

The derivative at $t=0$ being positive means
$$\RRe\langle Ay,A(a-y)\rangle-\mu^2\RRe\langle y,a-y\rangle>0,$$
or
$$\RRe\Big(\mu^2\langle y,a\rangle-\langle Aa,Ay\rangle\Big)<\mu^2\|y\|^2-\|Ay\|^2,$$
And one may exchange $a$ and $y$ here to find
$$\RRe\Big(\mu^2\langle y,a\rangle-\langle Aa,Ay\rangle\Big)<\mu^2\|a\|^2-\|Aa\|^2.$$
Replacing $y$ by $-y$, one gets similar inequalities for absolute values
\begin{eqnarray*}
&&\Big|\RRe\Big(\mu^2\langle y,a\rangle-\langle Aa,Ay\rangle\Big)\Big|<\mu^2\|y\|^2-\|Ay\|^2,\\
&&\Big|\RRe\Big(\mu^2\langle y,a\rangle-\langle Aa,Ay\rangle\Big)\Big|<\mu^2\|a\|^2-\|Aa\|^2.
\end{eqnarray*}
Which makes the LHS expression in (\ref{eq:ay}) positive, a contradiction.
\end{Prf}

\medskip

Let now $\mu>\mu'\ge0$, $I=[0,\mu]$,\,$I'=[0,\mu']\subset I$. Of course, every $I'$-subspace is \textit{ipso facto} an $I$-subspace, so is contained in some \textit{maximal} $I$-subspace.\pa

\begin{proposition}
Let $\mu>\mu'\ge0$, $I'=[0,\mu']$,\,\,$I=[0,\mu]$.

Consider a maximal $I'$-subspace $Y'$ contained in a maximal $I$-subspace $Y$. Then the quotient space $Y/Y'$ is isomorphic to $\mathcal{E}(I\setminus I')\mathcal{H}$.
\end{proposition}
\begin{Prf}
$\mathcal{E}(I)\mathcal{H}$ has the orthogonal decomposition $\mathcal{E}(I')\mathcal{H}\oplus\mathcal{E}(I\setminus I')\mathcal{H}$. The orthogonal projection here on $\mathcal{E}(I')\mathcal{H}$ will be the restriction to $\mathcal{E}(I')\mathcal{H}$ of the orthogonal projection from the whole space $\mathcal{H}$ on $\mathcal{E}(I')\mathcal{H}$.

Of course, $\mathcal{E}$ being a spectral decomposition, we have
$$\mathcal{E}(I')\mathcal{E}(I)=\mathcal{E}(I)\mathcal{E}(I')=\mathcal{E}(I').$$
We know that the restriction to $Y$ of the orthogonal projection $\mathcal{E}(I)$ from the whole space $\mathcal{H}$ on $\mathcal{E}(I)\mathcal{H}$ is an isomorphism with its image $\mathcal{E}(I)\mathcal{H}$. It maps the subspace $Y'$ to the subspace $\mathcal{E}(I)Y'\subset\mathcal{E}(I)\mathcal{H}$ there which the `further' orthogonal projection $\mathcal{E}(I')$ on $\mathcal{E}(I')\mathcal{H}$ will map isomorphically onto the latter.\pa

Hence, if we refer to $\mathcal{E}(I)\mathcal{H}$ as the product of $\mathcal{E}(I')\mathcal{H}$ and $\mathcal{E}(I\setminus I')\mathcal{H}$, then $\mathcal{E}(I)Y'$ is the `graph' of some bounded linear operator $B:\mathcal{E}(I')\mathcal{H}\to\mathcal{E}(I\setminus I')\mathcal{H}$.

This makes $Y/Y'$ isomorphic to the quotient space of $\mathcal{E}(I)\mathcal{H}$, identified with the product of $\mathcal{E}(I')\mathcal{H}$ and $\mathcal{E}(I\setminus I')\mathcal{H}$, over that `graph' $\{y, By\}$, and the latter is isomorphic to $\mathcal{E}(I\setminus I')\mathcal{H}$ via the map $(y,z)\mapsto z-By$.
\end{Prf}

\begin{corollary}\label{cor:max}
All maximal $[0,\mu]$-subspaces (with respect to $T$, equivalently to the positive semi-definite $A=\sqrt{T^\star T}$) are isomorphic, in particular are isomorphic to $\mathcal{E}([0,\mu])\mathcal{H}$ (hence have the same (Hilbert) dimension $\dim\mathcal{E}(I)$). And if $\mu\ge\mu'$ then all quotients of a maximal $[0,\mu]$-subspace over an included maximal $[0,\mu']$-subspace are isomorphic to $\mathcal{E}(]\mu',\mu])\mathcal{H}$.
\end{corollary}

\section{The Cardinal Number-Valued `Dimension Spectral Measure' as an Invariant}\label{s:card-meas}
As usual, by the \textit{(Hilbert) dimension} $\mathrm{dim}\,\mathcal{H}$ of a Hilbert space $\mathcal{H}$ we mean the cardinality of an orthonormal basis to $\mathcal{H}$ -- see \S\ref{s:dm-hil} (as is well known and easily proved, all such bases have the same cardinality).\pa

So, returning to our operator $T$ in $\mathcal{H}$ and the spectral decomposition $\mathcal{E}$ of $A=\sqrt{T^\star T}$, it defines a \textit{cardinal-number-valued ($\sigma$-additive) measure} on Borel subsets
$S\subset\mathbb{R}^+=[0,\infty[$, namely $\mathrm{dim}\,\mathcal{E}(S)\mathcal{H}$, which we call the
\textbf{dimension spectral measure} of $A$ and $T$.\pa

Clearly the dimension spectral measure determines a positive semi-definite operator $A$ up to replacing $A$ by $UAU^{-1}$, $U$ unitary.\pa

Thus Corollary \ref{cor:max} implies that \textit{for any $\mu\ge0$, all the dimensions of maximal $[0,\mu]$-subspaces are given by that dimension (cardinal number-valued) spectral measure $\mathrm{dim}\,\mathcal{E}(S)\mathcal{H}$ of $S=[0,\mu]$; and for any $\mu>\mu'\ge0$ all the dimensions of quotients of a maximal $[0,\mu]$-subspace over an included maximal $[0,\mu']$-subspace are given by the dimension spectral measure $\mathrm{dim}\,\mathcal{E}(S)\mathcal{H}$ of $S=]\mu',\mu]$}.\pa

But suppose we replace some bounded operator $T$ between Hilbert spaces by a topologically equivalent $T'$, i.e.\ $T'=V^{-1}TU$ where $U$ and $V$ are isomorphisms onto between Hilbert spaces. We may think about that as $T$ a continuous operator between locally convex linear topological spaces with different equivalent Hilbert norms defining the topology taken in them.\pa

Then there is some constant $K_1>0$ so that \textit{a maximal $[0,\mu]$-subspace $Y$ of $T$ with respect to one of the norms is a $[0,K_1\mu]$-subspace, hence included in a maximal $[0,K_1\mu]$-subspace, for some other}. And since by Corollary \ref{cor:max} the dimensions of quotients of such characterize the `dimension spectral measure' of the operators, we conclude that\pa

\begin{theorem}\label{th:dm-eq}
For operators between Hilbert spaces, to be topologically equivalent their dimension (cardinal number-valued) spectral measures must be equivalent in the following sense:\pa
\end{theorem}
\begin{definition}\label{def:dm-eq}
Two cardinal-number-valued ($\sigma$-additive) measures are called \textbf{equivalent} if there exists a constant $K\ge1$ such that:

The weight given to an interval $[0,\mu]$ (resp.\ $]\mu',\mu]$) by one of the measures is always $\le$ than the weight given to $[0,K\mu]$ (resp.\ $]\mu'/K,K\mu]$) by the other.
\end{definition}
And, of course, since linear operators, to be continuous must be bounded, the spectral measure, hence the dimension (cardinal number-valued) spectral measure of an operator (between Hilbert spaces) is supported on some bounded interval $[0,M]$.

\section{Focusing on Invariants to Characterize Topological Equivalence}
\renewcommand{\theenumi}{\roman{enumi}}
\begin{enumerate}
\item
Let $\mathcal{M}$ be some cardinal-number-valued ($\sigma$-additive) measure.\pa
%\begin{lemma}\label{le}

\item If a compact $\mathcal{K}\subset[0,\infty[$ has $\mathcal{M}(\mathcal{K})$ finite, while all its members $x\in \mathcal{K}$ are with $\mathcal{M}(\{x\})=0$ then $\mathcal{M}(\mathcal{K})=0$; By $\sigma$-additivity `compact' can be replaced here by `$\sigma$-compact', i.e.\ a countable union of compacta.
\begin{Prf}
We prove that every $x\in \mathcal{K}$ has an open neighborhood $]x-\eps,x+\eps[$ with $\mathrm{M}(\mathcal{K}\cap ]x-\eps,x+\eps[)=0$.

Indeed, $\mathrm{M}(\mathcal{K}\cap ]x-\eps,x+\eps[)$ is a decreasing function of $\eps$,\,$\eps>0$ with finite values, because $\mathrm{M}(\mathcal{K}\subset[0,\infty[)$ is finite. So it has a minimum value attained for some $\eps\le\eps_0$. But that implies that the weight of differences $\mathrm{M}(]x-\eps_0/(n-1),x-\eps_0/n]\cup[x+\eps_0/n,x+\eps_0/(n-1)[=0$, and by $\sigma$-additivity $\mathrm{M}(]x-\eps,x[\cup]x,x+\eps[=0$; Since by assumption $\mathrm{M}(\{x\})=0$, we have
$\mathrm{M}(]x-\eps,x+\eps[=0$.

Now, these neighborhoods for all $x\in \mathcal{K}$ are an open cover of $\mathcal{K}$; it has a finite subcover, whose union will have weight $0$; therefore $\mathcal{K}$ has weight $0$.
\end{Prf}

\item If a $\sigma$-compact $\mathcal{K}\subset[0,\infty[$ has $\mathrm{M}(\mathcal{K})$ finite, then there is a finite set $F\subset \mathcal{K}$ such that $\mathrm{M}(\mathcal{K}\setminus F)=0$.
\begin{Prf}
Let $F$ be the set of all $x\in \mathcal{K}$ with $\mathrm{M}(\{x\})>0$, the latter thus a positive integer. Clearly
$F$ has no more than $\mathrm{M}(F)\le\mathrm{M}(\mathcal{K})$ members, so it is finite. And then by the previous item $\mathcal{K}\setminus F$ must have weight $0$, being a $\sigma$-compact with every singleton of weight $0$.
\end{Prf}

\item If in a compact $\mathcal{K}\subset[0,\infty[$ every point has a neighborhood $U$ with $\mathrm{M}(\mathcal{K}\cup U)$ finite then $\mathcal{K}$ has finite weight.
\begin{Prf}
This is standard compactness argument: the neighborhoods with finite weight form an open cover of $\mathcal{K}$; it has a finite subcover; therefore the union of the latter has finite weight; hence so has $\mathcal{K}$.
\end{Prf}

\item Denote by $\mathbf{FIN}_{\mathcal{M}}$ the set of all $x\in[0,\infty[)$ with some neighborhood $]x-\eps,x+\eps[$ of finite weight, and by $\mathbf{INF}_{\mathcal{M}}$ its complement $\mathbf{INF}_{\mathcal{M}}:=[0,\infty[\setminus\mathbf{FIN}_{\mathcal{M}}$.

This makes $\mathbf{FIN}_{\mathcal{M}}$ the union of all open intervals of finite weight. Therefore it is an \textit{open} subset of
$[0,\infty[$, consequently a finite or countable union of open intervals there -- its connected components, and $\mathbf{INF}_{\mathcal{M}}$ is \textit{closed}.

By the above, any compact interval $[a,b]\subset\mathrm{FIN}_{\mathcal{M}}$ must have finite weight. Thus it has a finite subset where all the mass is concentrated. In any connected component (thus an open interval) $C$ of $\mathbf{FIN}_{\mathcal{M}}$, the mass will thus be concentrated in some discrete finite or infinite subset $\Sigma$, in which every point $x$ not the last (resp.\ first) will have an immediate successor (resp.\ predecessor) -- the biggest element of $\Sigma$ strictly less than $x$ (resp.\ the smallest element of $\Sigma$ strictly greater than $x$).

Therefore $\Sigma$ (of some component -- open interval -- of $\mathbf{FIN}_{\mathcal{M}}$) can be written as an increasing sequence, finite or infinite forward or backward or both, where if infinite it will converge to an extremity of $C$ (necessarily a member of $\mathbf{INF}_{\mathcal{M}}$).

\item
Returning to $\mathcal{D}$, the dimension spectral measure of an operator, the above $\Sigma$ will consist of
\textit{eigenvalues} of $A$, each with a finite-dimension space of eigenvectors, these spaces for different eigenvalues \textit{orthogonal}. And in such space consisting of eigenvector of some eigenvalue $\lambda$, say of dimension $n$, take an orthonormal basis of $n$ eigenvectors. Then $\Sigma$ (counting such $\lambda$ $n$ times) will turn into a \textit{non-decreasing sequence of eigenvalues} with a corresponding orthonormal set of eigenvectors.

\item
And whenever we have a sequence $(\lambda_k)$ of eigenvalues of $A$ and corresponding orthonormal eigenvectors $(v_k)$,\,\,$Av_k=\lambda_kv_k$, take some sequence of scalars $\gamma_k$ `uniformly log-bounded' (put otherwise, `bounded and bounded away from $0$'), i.e.\ there exists a $K>1$ such that for all $k$,
$$K^{-1}\le\gamma_k\le K.$$
Define an operator $B$ in the Hilbert space by $Bv_k:=\gamma_kv_k$ and $B$ the identity on the orthogonal complement of the subspace spanned by the $v_k$.

Then $B$ is invertible -- its inverse obtains as the analog replacing each $\gamma_k$ by $\gamma_k^{-1}$. $B$ also commutes with $A$, and their product (thus, topologically equivalent with $A$) is an operator with each $v_k$ now an eigenvector of $\gamma_k\lambda_k$ -- in the dimension spectral measure the $1$-mass was `transferred' from $\lambda_k$ to $\gamma_k\lambda_k$. Hence

\begin{corollary}
Suppose $\gamma_k$ is `uniformly log-bounded' as above, then two operators $T$ and $T'$ with their dimension (cardinal-number-valued) spectral measures obtained one from the other by `transferring $1$-masss from $\lambda_k$ to $\gamma_k\lambda_k$', must necessarily be topologically equivalent.
\end{corollary}

We will sometimes refer to such \textit{transfer of masses, in a cardinal-number-valued ($\sigma$-additive) measure, to uniformly log-bounded-ly distant places} as a \textbf{uniformly log-bounded shift} in the measure.

\item
Now let us be given some countable partition of $]0,\infty[$ into Borel subsets:
$$]0,\infty[=\cup_k S_k,\quad\mathrm{the}\,S_k\,\mathrm{disjoint}\,\mathrm{Borel}\,\mathrm{sets},$$
and some corresponding sequence of points $a_k\in]0,\infty[$, such that \textit{the log-diameters of $S_k\cup\{a_k\}$ are uniformly bounded}: there is a $K>1$ such that for all $k$,
$$y,z\in S_k\cap\{x_k\}\Rightarrow K^{-1}\le y/z \le K.$$
As an example, let $\beta>1$ and take $S_k=]b\beta^k,b\beta^{k-1}],\,\,k\in\mathbb{Z}$ (or $S_k=[b\beta^k,b\beta^{k-1}[,\,\,k\in\mathbb{Z}$) and $a_k=a\beta^k,
\,\,k\in\mathbb{Z}$ for some fixed $a,b>0$.\pa

Returning to our operator $T$ between Hilbert spaces and (the positive semi-definite) $A:=\sqrt{T^\star T}$, with spectral decomposition $\mathcal{E}$. Recall that
$$A=\int_{]0,+\infty[}\lambda\,d\mathcal{E}(\lambda).$$
Let the scalar function $\varphi$ on $]0,\infty[$ be defined by $\varphi(\lambda):=a_k/\lambda$ when $\lambda\in S_k$.
Then by our assumptions we will have $K^{-1}\le\varphi(\lambda)\le K$ for $\lambda\in]0,\infty[$.

Therefore if we define an operator from $\mathcal{H}$ to itself:
$$B:=\mathcal{E}(\{0\})+\int_{]0,+\infty[}\varphi(\lambda)\,d\mathcal{E}(\lambda),$$
then $B$ is invertible -- its inverse obtains as the analogous integral replacing the integrand by its inverse. $B$ also commutes with $A$, and their product (thus, topologically equivalent with $A$) is
$$BA:=\int_{]0,+\infty[}\lambda\varphi(\lambda)\,d\mathcal{E}(\lambda)=
\sum_{k}\int_{S_k}a_k\,d\mathcal{E}(\lambda)=
\sum_{k}a_k\mathcal{E}(S_k).$$
In particular, we found that the operators $T$ and $A$ are topologically equivalent to \textit{a (positive semi-definite) operator whose spectral decomposition is supported on $0$ and the $a_k$, and moreover depends only on, in fact is a combination of, the $a_k$'s and the value $\mathcal{E}$ gives to $(\{0\})$ and to the $S_k$'s.}\pa

Therefore\pa

\begin{corollary}
With $S_k$ and $a_k$ as above, two operators $T$ and $T'$ with their dimension (cardinal-number-valued) spectral measure giving the same respective weights to $(\{0\})$ and to the $S_k$'s must necessarily be topologically equivalent.
\end{corollary}

Or, to put it another way,\pa

\textit{For an operator between Hilbert spaces, up to replacing the operator by a topologically equivalent one, we may, in the dimension spectral measure, `transfer' all the mass of $S_k$ into $a_k$ for each $k$}.\pa

Which is another case of a \textit{uniformly log-bounded shift} --
`transfer of mass (in a cardinal-number-valued ($\sigma$-additive) measure) to uniformly log-bounded-ly distant places'.
\end{enumerate}

\section{Towards a Classification, `Eigenvalue Sequences with Infinity Intermissions'}
We are given our operator $T$ between Hilbert spaces, topologically equivalent to the positive semi-definite $A:=T^\star T$, with its dimension spectral measure $\mathcal{D}$ -- the latter characterizing it up to  isometry, i.e.\ conjugation by a unitary, hence \textit{a fortiori} up to topological equivalence.\pa

For $\mathcal{D}$ we will have the open (subset of $]0,\infty[$) $\mathbf{FIN}_{\mathcal{D}}$, thus an at most countable union of open intervals -- its connected components; and its complement -- the closed (in $]0,\infty[[$ $\mathbf{INF}_{\mathcal{D}}=]0.\infty[\setminus\mathbf{FIN}_{\mathcal{D}}$.\pa

Then, if necessary making uniformly log-bounded shifts in the dimension spectral measure $\mathcal{D}$, thus replacing an operator by a topologically equivalent, one may
\begin{itemize}
\item Fixing a $\beta>1$, `transfer' all masses to, say, $0$ and $\beta^k,\,k\in\mathbb{Z}$.
\item Recall that the mass is necessarily supported in some interval $[0,a]$ since the operator is bounded. And we can make $a$ an arbitrarily smaller $a'<a$ by `transferring' the mass of $[a,a'[$ into $a'$.
\item Recall the decomposition of $]0,\infty[$ into $\mathbf{FIN}_{\mathcal{D}}$ and $\mathbf{INF}_{\mathcal{D}}$.
    In $\mathbf{INF}_{\mathcal{D}}$ no point $x$ has a neighborhood with finite weight. Thus either $]x-\eps,x]$ has infinite weight for every $\eps>0$, or $[x,x+\eps[$ are so (or both). Then, applying a uniformly log-bounded shift, we may assume $\{x\}$ itself has infinite weight.
\item And if in a component $C$ of $\mathbf{FIN}_{\mathcal{D}}$, the set $\Sigma$ of eigenvalues turns to be infinite forward or backward, thus tending to an extremity $a$ of $C$, necessarily $a\in\mathbf{INF}_{\mathcal{D}}$, the infinite part from any fixed place can be `transferred' to $a$, thus `truncating' the sequence to leave only a finite sequence of eigenvalues, and surely making $\mathcal{D}(\{a\}$ infinite.
\item More generally, a point $a$ with \textit{infinite weight} may serve to `erase', by uniformly log-bounded shift (`transferring them to $a$), any discrete points with finite weight uniformly log-bounded distanced from $a$.
\end{itemize}

So, after some uniformly log-bounded shift, one may assume the mass concentrated in $0,\lambda_k$, the $\lambda_k$ a (given in advance) strictly decreasing sequence, tending to $0$ and \textit{with log-distances uniformly bounded below}. i.e.\ $\lambda_k/\lambda_{k+1}\ge K$ for some fixed $K>1$ -- as a particular case we may choose
$\lambda_k:=\beta^{k_0-k}$, $\beta>1$ given.\pa

Some points $0,\lambda_k$ will have an infinite-cardinal weight.\pa

The other, with finite weight, may be arranges as a \textit{non-increasing} sequence of eigenvalues, counting $n$ times an eigenvalue $\lambda$ with multiplicity $n$, i.e.\ with weight $n$ in the dimension spectral measure (and of course $\lambda$'s with weight $0$ not appearing). Recall that this non-increasing sequence is such that quotients of a term by its successor, if not $1$, are \textit{uniformly bounded away from $1$}.\pa

We call that an \textbf{eigenvalue sequence with infinity intermissions}.\pa

To repeat, up to topological equivalence, making a uniformly log-bounded shift if necessary, every operator is characterized by such eigenvalue sequence with infinity intermissions.\pa

Also, by uniformly log-bounded shift, one can `erase' some or all eigenvalues in the eigenvalue sequence with uniformly bounded log-distance from points $0,\lambda_k$ with infinite weight.\pa

Recall that Theorem \ref{th:dm-eq} with Definition \ref{def:dm-eq} gave a \textit{necessary condition} for operators between Hilbert spaces, with corresponding dimension spectral measures $\mathcal{D}$ and $\mathcal{D}'$, to be
topologically equivalent: namely, that there exists a constant $K\ge1$ such that:\pa

The weight given to an interval $[0,\mu]$ (resp.\ $]\mu',\mu]$) by one of the measures is always $\le$ than the weight given to $[0,K\mu]$ (resp.\ $]\mu'/K,K\mu]$) by the other.\pa

What does that say in our scenario?\pa

If $\lambda>\lambda'$ are terms in the eigenvalue sequence with \textit{no infinity intermissions between them} and with $\lambda'$\,\,$n$ steps after $\lambda$, then clearly for any $\eps>0$,\pa

The weight of $]\lambda'-\eps,\lambda]\ge n\ge$ the weight of $]\lambda',\lambda-\eps]$.\pa

Therefore the necessary condition translates here to:\pa

In the eigenvalue sequences of the operators, starting from some $\lambda$, even allowing a uniformly log-bounded shift in it -- and then going $n$ steps ($n$ finite) in the sequence for either of the operators without encountering infinity intermissions -- will give for one operator results at uniformly log-bounded distance (that is, quotient) from the other operator.\pa

And also, that points with infinite cardinal weight for one of the operators must have, at uniformly log-bounded distance, points with no less weight for the other operator.\pa

Where, recall, as long as we are in uniformly log-bounded distance from points with infinite weight, eigenvalues in the sequence may be erased (by uniformly log-bounded shift).\pa

Thus in starting, as above, from some $\lambda$, even allowing a uniformly log-bounded shift in it, then going $n$ steps ($n$ finite) in the sequence -- thing matter only as long as we are not at uniform log-bounded distance from
points with infinite weight, for one operator consequently for the other -- our necessary condition then says that going $n$ steps for one operator or the other will keep us at uniformly log-bounded distance.\pa

And clearly if all that is satisfied, the eigenvalue sequence with infinity intermissions for one operator obtains from that for the other by some uniformly log-bounded shift. Therefore\pa

\begin{theorem}\label{th:cond}
For two operators between Hilbert spaces, a necessary and sufficient condition to be topologically equivalent is that their `eigenvalue sequence with infinity intermissions' can be obtained one from the other by uniformly log-bounded shift.
\end{theorem}

\textbf{discussion}
\begin{itemize}
\item Of course, $\mathcal{D}(\{0\})$ -- the (Hilbert) dimension of the kernel, is an invariant for topological equivalence.
\item So are $\mathcal{D}([0,\infty[)$ -- the (Hilbert) dimension of the space and $\mathcal{D}(]0,\infty[)$ -- the (Hilbert) dimension of the image, in particular when they are finite.
\item When \textit{the number of eigenvalues between the infinity intermissions is uniformly bounded} we may, by a uniform log-bounded shift assume that \textit{all points with non-zero weight have infinite weight}.
\item The case where \textit{there are no infinity intermissions} is clearly the case of \textit{a compact operator} -- obviously a topological equivalence invariant.
\item What may yield non-topologically-equivalent operators is, for instance, unbounded `chunks' in the eigenvalue sequence between `infinities'-- recall these are non-increasing, with quotients of a term by its successor, if not $1$, \textit{uniformly bounded away from $1$}. Hence what can differentiate here are thing like \textit{unbounded lengths of rows of equal terms, i.e.\ eigenvalues with multiplicity}; and \textit{kog-unbounded `jumps' from an eigenvalue term to the next}.
\item Also, especially when the Hilbert space is non-separable, \textit{which infinite cardinalities appear and their distribution in the infinity intermissions} -- again inasmuch as not `uniformly bounded',
\end{itemize}

\appendix

\section{Dimension in Hilbert Spaces}\label{s:dm-hil}
Recall that, as usual, by the \textit{dimension} $\mathrm{dim}\,\mathcal{H}$ of a Hilbert space $\mathcal{H}$ we mean the cardinality of an orthogonal basis to $\mathcal{H}$ (as is well known and easily proved, all such bases have the same cardinality).\pa

\begin{proposition} \label{pr:HilDim}
If $\mathcal{H}$ and $\mathcal{H}'$ are Hilbert spaces and $L:\mathcal{H}\to\mathcal{H}'$ is an injective bounded linear operator, then $\dim\mathcal{H}\le\dim\mathcal{H}'$.
\end{proposition}
\begin{Prf}
Let $\delta=\dim\mathcal{H}$ (resp.\ $\delta'=\dim\mathcal{H}'$) and let $(e_i)_{i\in\delta}$ be an orthonormal basis
to $\mathcal{H}$.\pa

If $\delta$ is finite then $L\mathcal{H}\subset\mathcal{H}'$ is also a finite-dimensional space of dimension $\delta$
hence $\delta\le\dim\mathcal{H}'=\delta'$.\pa

Assume now that $\delta$, hence $\delta'$, is infinite.\pa

$L$ is strongly, hence weakly, continuous, and is injective. Therefore it defines a homeomorphism, with respect to the
weak topologies, from the unit ball of $\mathcal{H}$, which is weakly compact, to its image, consequently also weakly
compact.\pa

The set $\{e_i\}$ has weakly the discrete topology in the unit ball of $\mathcal{H}$, therefore also the set $\{Le_i\}$
has weakly the discrete topology in $\mathcal{H}'$.\pa

Now, $\mathcal{H}'$ has an open basis to the weak topology of cardinality $\max(\delta',\aleph_0)=\delta'$ (consisting of
weakly open sets defined by inequalities on the inner product with finite rational combinations of the orthonormal basis).
Since $\{Le_i\}$ is weakly discrete, every $Le_i$ belongs to some member $U_i$ of that basis containing it and no other
$Le_j$,\,$j\ne i$. Of course, all the $U_i$'s must be different. Thus $\delta$, the cardinality of $\{Le_i\}$, is less
than $\delta'$.
\end{Prf}

\begin{corollary}\label{cr:HilDim}
Two isomorphic Hilbert spaces (i.e.\ such that there is a 1-1 and onto, i.e.\ invertible, operator between them), have the same Hilbert dimension, hence are, in fact, isometric.
\end{corollary}

\section{In Hilbert Space, the Subspace of Vectors Where the Norm is Attained}
We are in the scenario of \S\ref{s:A} and \S\ref{s:iso}.\pa

For a vector $x\in\mathcal{H}$, we have
\begin{equation}\label{eq:norm}
\|Tx\|^2=\langle A^2x,x\rangle=\left\langle\int_{[0,+\infty[}\lambda^2\,d\mathcal{E}(\lambda)x\,\,,\,\,x\right\rangle=
\int_{[0,+\infty[}\lambda^2\,\left\langle\,d\mathcal{E}(\lambda)x,x\right\rangle.
\end{equation}
This means that $\|T\|$ is `the upper end of the support of the spectral measure $\mathcal{E}$'. Indeed, $\mathcal{E}(]\|T\|,\infty[)=0$, since otherwise the subspace $\mathcal{E}(]\|T\|,\infty[)\mathcal{H}$ would be a non-zero subspace, thus including some $x\ne 0$ for which by (\ref{eq:norm}) $\|Tx\|$ would be $>\|T\|\|x\|$; while on the other hand for any $\eps>0$ $\mathcal{E}(]\|T\|-\eps,\infty[)\ne0$, otherwise by (\ref{eq:norm}) for every $x\in\mathcal{H}$\,\,$\|Tx\|\le(\|T\|-\eps)\|x\|$, contrary to the definition of $\|T\|$.\pa

Moreover, one finds that for a vector $x\in\mathcal{H}$ to satisfy $\|Tx\|=\|T\|\|x\|$, one must have, by (\ref{eq:norm})
\begin{eqnarray*}
&&\int_{[0,\|T\|]}\lambda^2\,\left\langle\,d\mathcal{E}(\lambda)x,x\right\rangle=
\|T\|^2\int_{[0,+\infty[}\left\langle\,d\mathcal{E}(\lambda)x,x\right\rangle,\\
&&\int_{[0,\|T\|[}(\|T\|^2-\lambda^2)\,\left\langle\,d\,\mathcal{E}(\lambda)x,x\right\rangle=0,
\end{eqnarray*}
which, since $\|T\|^2-\lambda^2$ is \textit{strictly positive} in $[0,+\|T\|[$, can happen only if
$\|\mathcal{E}([0,+\|T\|[)x\|^2=\int_{[0,+\|T\|[}\left\langle\,d\mathcal{E}(\lambda)x,x\right\rangle=0$, meaning that
$x\in\mathcal{E}(\{\|T\|\})\mathcal{H}$. Thus we reach

\begin{Cor}\label{cor:norm}
For any operator $T$ on a Hilbert space $\mathcal{H}$, and $\mathcal{E}$ the spectral measure of
$A:=\sqrt{T^\star T}$ as above, the set of vectors $x\in\mathcal{H}$ satisfying $\|Tx\|=\|T\|\|x\|$ is the \textit{(closed) subspace} $\mathcal{E}(\{\|T\|\})\mathcal{H}$.
\end{Cor}

\end{document}